\documentclass[12pt]{article}
\usepackage{amsmath,amsfonts,amssymb}

\def\diag #1{{\rm diag}\left(#1\right)}
\def\1{\mbox{\bf 1}}

\def\r{\rightarrow}
\def\R{\mathbb{R}}

\def\iy{\infty}
\def\qed{\hfill\vrule height5pt width5pt depth0pt}
\def\one #1{1_{\{#1\}}}

\def\diag #1{{\rm diag}\left(#1\right)}
\def\eq #1{(\ref{eq:#1})}
\def\Th #1{Theorem~\ref{Th:#1}}

\def\Ex #1{Example~\ref{Ex:#1}}

\newcommand{\bea}{\begin{eqnarray}}
\newcommand{\eea}{\end{eqnarray}}
\newcommand{\Bea}{\begin{eqnarray*}}
\newcommand{\Eea}{\end{eqnarray*}}
\newcommand{\proof}{\noindent {\bf Proof:\ }}

\newtheorem{Definition}{Definition}[section]
\newtheorem{Theorem}[Definition]{Theorem}
\newtheorem{Lemma}[Definition]{Lemma}
\newtheorem{Proposition}[Definition]{Proposition}
\newtheorem{Corollary}[Definition]{Corollary}
\newtheorem{Remark}[Definition]{Remark}
\newtheorem{Example}[Definition]{Example}

\begin{document}

\title{Asymptotic irrelevance of initial conditions for Skorohod reflection mapping on the nonnegative orthant}

\author{
Offer Kella\thanks{Department of Statistics; The Hebrew University
of Jerusalem; Jerusalem 91905; Israel ({\tt
Offer.Kella@huji.ac.il}).} \thanks{Supported in part by grant No.
434/09 from the Israel Science Foundation and the Vigevani Chair in
Statistics.}
 and S. Ramasubramanian\thanks{Theoretical Statistics and Mathematics Unit; Indian Statistical Institute; 8th mile, Mysore Road; Bangalore - 560 059; India ({\tt ram@isibang.ac.in}).}
}
\date{Dec. 20, 2011}
\maketitle
\begin{abstract}
A reflection map, induced by the deterministic Skorohod problem on
the nonnegative orthant, is applied to an $\mathbb{R}^n$ valued
function $X$ on $[0,\iy)$ and then to $a+X$, where $a$ is a
nonnegative constant vector. An question that has been posed over
15 years ago is under what conditions the difference between the two
resulting regulated functions converges to zero for any choice of
$a$ as time diverges. This in turn implies that if one imposes
enough stochastic structure that ensures that the reflection map
applied to a multidimensional process $X$ converges in distribution,
then it will also converge in distribution when it is applied to
$\eta+X$ where $\eta$ is any almost surely finite valued random
vector that may even depend on the process $X$. In this paper we
obtain a useful equivalent characterization of this property. As a
result we are able to identify a natural sufficient condition in
terms of the given data $X$ and the constant routing matrix. A
similar necessary condition is also indicated. A particular
implication of our analysis is that under additional stochastic assumptions,
asymptotic irrelevance of the initial
condition does not require the existence of a stationary distribution.
As immediate corollaries of our (and earlier) results we conclude that under the
natural stability conditions, a reflected L\'evy process as well as Markov additive process
has a unique stationary distribution and converges in distribution to this stationary distribution for
every initial condition.
Extensions of the
sufficient condition are then developed for reflection maps with
drift and routing coefficients that may be time and state dependent;
some implications to multidimensional insurance models are briefly
discussed.


\bigskip\noindent {\bf AMS Subject Classification (MSC2010):} Primary 90B15; Secondary 60G17, 60G51,
60K99, 90B10, 91B30.

\end{abstract}

\section{ Introduction and background}
\setcounter{equation}{0}

Let $X=\{X(t)|\ t\ge 0\}$ be an $\mathbb{R}^n$ valued c\`adl\`ag
(right continuous left limit) function, $R=I-P^t$, where $^t$
denotes transposition and $P$ is nonnegative with spectral radius
less than one, that is, equivalently with $P^n\r0$ as $n\r\infty$.
The matrix $R$ is called an $M$-{\em matrix} (e.g., \cite{BP1994}).
Consider
\begin{equation}W(X)(t)=X(t)+RL(X)(t)\end{equation}
where $(W,L)$ is the unique
(reflection) mapping of $X$ satisfying for each $i=1,\ldots,n$:
\begin{description}
\item{S1} $L_i(X)$ is nondecreasing, nonnegative and right continuous.
\item{S2} $W_i(X)$ is nonnegative.
\item{S3} $\int_{[0,\infty)}\one{W_i(X)(t)>0}dL_i(X)(t)=0$.
\end{description}
We note that this is also often referred to as a Skorohod problem in
$[0,\infty)^n$ and from time to time this is what we will call it,
especially in the setting of (\ref{eq:skorohod}) below. It is
mentioned in passing that in \cite{K2006} an alternative but
equivalent and somewhat more microscopic definition of the
reflection map is given and it is shown that it is well defined
under weaker conditions on $X$. That is, the reflection map is well
defined even if $X$ is not necessarily c\`adl\`ag, but rather lower
semicontinuous from the right.

Now, with the current notations, Theorem~6 of \cite{KW1996} can be rephrased as follows.

\begin{Theorem}\label{Theorem:1.1}
Assume that $X^1,X^2$ are c\`adl\`ag, $X^2-X^1$ is nonnegative and nondecreasing and denote
$X^3=X^1+RL\left(X^2\right)$. Let $e=(1,\ldots,1)^t$. Then

\begin{description}

\item{\rm(i)} $W\left(X^2\right)(t)\ge W\left(X^1\right)(t)$ for
all $t\ge 0$.

\item{\rm(ii)} $L\left(X^3\right)=L\left(X^1\right)-L\left(X^2\right)$
(in particular nonnegative and nondecreasing) and is dominated above by
$R^{-1}\left(X^2-X^1\right)$.

\item{\rm(iii)}
$R^{-1}\left(W\left(X^2\right)-W\left(X^1\right)\right)$
(resp. $e^t\left(W\left(X^2\right)-W\left(X^1\right)\right)$) is
nonnegative and dominated above by $R^{-1}\left(X^2-X^1\right)$ (resp.
$e^t\left(X^2-X^1\right)$).

\item{\rm(iv)} When $X^2(t)=a+X^1(t)$, where $a$ is a nonnegative constant
vector, then
$R^{-1}\left(W\left(X^2\right)-W\left(X^1\right)\right)$,
hence $e^t\left(W\left(X^2\right)-W\left(X^1\right)\right)$, is
nonincreasing.

\end{description}
\end{Theorem}

 It is mentioned that in \cite{KW1996} $P$ was also assumed to be
 substochastic.
 However, it can easily be checked that this is not needed in the
 proofs hence Theorem~\ref{Theorem:1.1} is valid without this additional restriction. In fact, in a
 much more general setting, similar results were obtained a few years
 later in Theorem~4.1 of \cite{R2000}. In particular, the model
 considered there is of the form

\begin{eqnarray}\label{eq:skorohod}
W(X)(t)\nonumber&=&X(t)+\int_0^t
b\left(u,L(X)(u-),W(X)(u-)\right)du\\
\\&&+\int_{[0,t]}R\left(u,L(X)(u-),W(X)(u-)\right)dL(X)(u)\
.\nonumber
\end{eqnarray}
where $b$ and $R$ satisfy certain regularity conditions.

These results can be used, applying a Loynes' type argument, to
establish (see \cite{KW1996,K1997}) that when $X$ is a c\`adl\`ag
stochastic process having stationary, ergodic (but not necessarily
independent!) finite mean increments and $R^{-1}E(X(1)-X(0))<0$,
then $W(X)$ is stable in the sense that it is tight and it has a
stationary version. That is, there is a random vector $\xi$ such
that $W(\xi+X)$ is a stationary process. For the stationary version
$L(\xi+X)$ has stationary increments. Also, if in addition $X(0)=0$,
then $W(X)(t)$ is stochastically increasing in $t$, hence has a
limiting distribution which coincides with the distribution of $W(\xi+X)(s)$ for each $s$.
What is missing from this theory is what happens when $X(0)$ is not
zero. Does there exist a limiting distribution for any initial
$X(0)$? Does it depend on $X(0)$? Is the stationary version unique? These questions are very relevant
especially in the context queueing networks; see \cite{KW1996,K1997,KL2004}
and references therein.

In Section 2 of this paper we establish necessary and sufficient conditions for
\begin{equation}W(a+X)(t)-W(X)(t)\r0\end{equation} as $t\r\infty$
for every nonnegative vector $a$.
 This problem was open \cite{KW1996,KL2004} for a
while and its resolution also resolves the issues raised at the
end of the last paragraph. Special cases are multidimensional L\'evy processes
and multidimensional Markov additive processes with countable (or finite) state space modulation.
See Theorem~\ref{Th:5}, Corollary~\ref{Cor:Levy} and Corollary~\ref{Cor:MAP} in what follows.

In Section 3 we explore possible extensions for a more general set up.
Regulated/ reflected processes constrained to stay in the nonnegative orthant
with time and general space dependent coefficients (that is, with drift
vector $b(\cdot)$ and reflection matrix $R(\cdot)$ depending on the variables
$(t,y,z)$ with the notations of Section 3) have been studied
by various authors. In particular in the context of queueing networks and 
multidimensional insurance models. If the dispersion
is constant, but the drift and reflection are not necessarily
constant, solution to the stochastic problem can be obtained by
solving the associated deterministic Skorohod problem path-by-path;
so \eq{skorohod} has also received attention; see
\cite{DI1993,MP1998,R2000,R2011a} and references therein. In
connection with insurance models, it has been argued in
\cite{R2006,R2011a} that non-constant drift and reflection arise in
a natural fashion.

In one dimensional actuarial risk theory, the ruin problem has a central
role. Under certain natural conditions, it is known that the asymptotic
nature of the ruin is the same or similar independent of the initial capital (see
\cite{EK1997,RS1999}). It has been shown in \cite{R2011b} that the
regulated process in the orthant hitting state $0$ is the
appropriate notion of ruin for the multidimensional set up
considered in \cite{R2006,R2011a}. Extensions considered in Section
3 at the level of sample paths might hopefully lead to a better
understanding of the multidimensional ruin problem; see \Ex{Re1}.

We may now indicate how our work stands in relation to earlier
studies. Asymptotic properties of constrained non-degenerate
diffusion processes or jump-diffusion processes have been studied by
many authors; see \cite{AB2002} for a fairly comprehensive list. As
these processes are strong Markov and their infinitesimal generators
are non-degenerate differential or integro-differential operators,
appropriate Lyapunov functions, when they exist, can be used to establish existence of a unique stationary distribution
independent of the initial state or distribution. In our set up,
such versatile tools are not available in general, as there may not
even be an underlying probabilistic framework. Our main concern is the
asymptotic pathwise irrelevance of the initial condition. As will be seen below,
if in addition a stochastic structure is imposed that guarantees existence
of a stationary distribution, then uniqueness also follows. Insensitivity to the
initial condition, however, can happen even when there is no stationary distribution;
see \Ex{BM}.

During the review process our attention has been drawn to
\cite{BC2002,BC2006}, where a problem similar to the one in Section
2 has been studied for reflected Brownian motion in a bounded planar
domain with normal reflection at the boundary; see also
\cite{CL1990} for an earlier work in that direction. Let
$U(\cdot),\hat{U}(\cdot)$ be reflected Brownian motions in a bounded
domain $D\subset \R^2$ with normal reflection starting respectively
at $x,\hat{x}\in D.$ Asymptotic behaviour of $|U(t)-\hat{U}(t)|$ as
$t\r\iy$ has been investigated in terms of nice geometric/ analytic
properties related to the domain. The infinitesimal generator in
this case is the two dimensional Laplacian with Neumann boundary
condition and the normalized Lebesgue measure is the stationary
distribution. Tools from spectral theory of self-adjoint operators
are used and beautiful connections with the geometry of $D$ are
elucidated in \cite{BC2002,BC2006,CL1990}. On the other hand, our investigation at the sample path level
concerns oblique reflection as well, though restricted to a
multidimensional orthant which is unbounded. Also, as mentioned
earlier, our focus is on insensitivity to the initial
state and that too at the level of sample paths. There need not
even be any underlying probabilistic framework, let alone a
stationary or a Markovian structure or associated mathematical
machinery. As a consequence of our sample path results applied to
stochastic processes, it turns out that insensitivity to the  initial state can
happen with probability one even in cases where the reflected
process may be null recurrent or transient; see \Ex{BM}. A similar
almost sure result is given in \Ex{Renewal} for a non-Markovian
process.

\section{Main results}
\setcounter{equation}{0}

It is reemphasized that due to Theorem 1.1,  $R^{-1}(W(a+X)-W(X))$
is nonincreasing and $W(a+X)-W(X)$ is nonnegative and thus, as $R$
has an inverse, namely $R^{-1}=\left(\sum_{i=0}^\infty
P^n\right)^t$, the limit of $W(a+X)-W(X)$ necessarily exists and is
finite. The question is whether it must be zero or not.

For the case where the states $1,\ldots,n$ can be reordered in such
a way that $p_{ij}=0$ for $i>j$, which is called a {\em feedforward}
structure, the problem is relatively easy to solve. In this case,
under appropriate natural conditions, there actually exists a finite
coupling time from which the process that starts from zero and the
one that starts from $a$ become identical. See Theorem~4.1 of
\cite{K1997} for details. Special cases are, of course, the one
dimensional, tandem ($p_{i,i+1}=1$ for $i=1,\ldots,n-1$)  and
parallel ($P=0$) cases, which have been given special treatments in
the literature, but are of no relevance to solving the general
problem, hence references are omitted.

Before we continue, we note that without loss of generality it may
be assumed that $p_{ii}=0$. The reason is that with
$D=\diag {p_{11},\ldots,p_{nn}}$,
\begin{equation}
\tilde W(X)=\tilde X+\tilde R \tilde L(X)
\end{equation}
where $\tilde Y=(I-D)^{-1/2}Y$ for $Y=X,W,L$ and $\tilde R=I-\tilde
P$ with \begin{equation}\tilde
P=(I-D)^{-1/2}(P-D)(I-D)^{1/2}\end{equation} is nonnegative and has
the same eigenstructure as $P-D$ which is also nonnegative and also
has a spectral radius less than one, as it is bounded above by $P$.
Note that since $D^n$ is bounded above by the diagonal of $P^n$ then
necessarily $p_{ii}<1$ for all $i$.

The following is the main result of this paper.
\begin{Theorem} \label{Theorem:2.1}
For any given c\`adl\`ag $X$ and every $i$, the following conditions
are equivalent:
\begin{description}
\item{\rm(i)} $\lim_{t\r\infty}[W_i(a+X)(t)-W_i(X)(t)]= 0$ for all $a\ge 0$.
\item{\rm(ii)} $\lim_{t\r\infty}L_i(a+X)(t)=\infty$ for  some $a\ge 0$.
\item{\rm(iii)} $\lim_{t\r\infty}L_i(a+X)(t)=\infty$ for all $a\ge 0$.
\end{description}
\end{Theorem}

\newpage
\proof
\begin{description}
\item{(ii)$\Leftarrow$(iii)} Obvious.
\item{(ii)$\Rightarrow$(iii)} Since
\begin{equation}R^{-1}W(a+X)-R^{-1}W(X)=R^{-1}a+L(a+X)-L(X)\ ,\end{equation}
by (ii) and (iv) of Theorem 1.1 note that
\begin{equation} L(a+X)\le L(X)\le L(a+X)+R^{-1}a\ ,\end{equation}
for all $a\geq 0.$ Now, if (ii) holds for some $a\geq 0,$ then by
the left inequality above it also holds for $a=0$ and by the right
inequality it holds for all $a\ge0$.
\item{(i)$\Leftarrow$(iii)} From (i) and (iv) of Theorem 1.1 we have that $W(a+X)\ge W(X)$ and, as discussed earlier, that \begin{equation}\lim_{t\r\infty}(W(a+X)(t)-W(X)(t)) \end{equation}
    exists. We denote this limit by $\zeta^a$ and observe that it must be
nonnegative. Assume that (i) does not hold, so that $\zeta^a_i>0$
for some $a$. Then there is some $T$ such that for $t\ge T$ we have
that
\begin{equation}
W_i(a+X)(t)-W_i(X)(t)> \frac{\zeta_i^a}{2}
\end{equation}
and since $W_i(X)(t)\ge0$ for all $t\ge0$ then in fact
$W_i(a+X)(t)>\zeta_i^a/2$ for all $t\ge T$. Since $\int_0^\infty
\one{W_i(a+X)(t)>0}dL_i(a+X)(t)=0$ it follows that
$L_i(a+X)(t)=L_i(a+X)(T)$ for all $t\ge T$. This contradicts (iii)
and thus (iii) implies (i).
\item{(i)$\Rightarrow$(iii)} Assume that (iii) does not hold and that, for some $a\ge 0$, $L_i(a+X)$ is bounded. Then from (ii)$\Leftrightarrow$(iii) it follows that $L_i(X)$ is bounded, say, by $b\ge0$. With $e_i$ being the unit vector with $1$ in the $i$th coordinate and zero elsewhere, we have from (ii) of Theorem 1.1 that $L(be_i+X)\le L(X)$ which implies that
\begin{eqnarray}\label{eq:nonnegative}
0\le W_i(X)&=&X_i+L_i(X)-\sum_jp_{ji}L_j(X)\nonumber\\ &\le &X_i+b-\sum_jp_{ji}L_j(be_i+X)\\
&=&(be_i+X-P^tL(be_i+X))_i\nonumber
\end{eqnarray}
Recall that for c\`adl\`ag $Y$, $L(Y)$ is the unique solution of the equations
\begin{equation}
L_k(Y)(t)=-\inf_{0\le s\le
t}\left[\left(Y_k(s)-\sum_jp_{jk}L_j(Y)(s)\right)\wedge 0\right],\begin{array}{c}1\leq k\leq n,\\ t\geq 0.\end{array}
\end{equation}
Set $Y=be_i+X,$ use the fact that $p_{ii}=0$ and the inequality in (\ref{eq:nonnegative}) to conclude that
$L_i(be_i+X)\equiv 0$. For all $c> b$ we also have from (ii) of
Theorem 1.1 that $L(ce_i+X)\leq L(be_i+X)$ and in particular the
$i$th coordinate is zero with both $c$ and $b$. Thus,
\begin{equation}\begin{array}{l}
W_i(ce_i+X)-W_i(be_i+X)\nonumber\\
\qquad=c-b-\sum_jp_{ji}(L_j(ce_i+X)-L_j(be_i+X))\\ \qquad\ge
\nonumber c-b>0\ .
\end{array}
\end{equation}
Therefore, subtracting and adding $W_i(X)$ and letting $t\r\infty$
we have that $\zeta^{ce_i}_i-\zeta^{be_i}_i\ge c-b>0$ and (recall $\zeta^{be_i}\ge
0$) in particular $\zeta^{ce_i}_i>0$ which contradicts (i). \qed
\end{description}

\begin{Remark}
\rm We note that we were not careful about assuming that $X(0)=0$ or even that $X(0)\ge 0$. This is because
\[L(X)(t)-L(X)(0)=L(W(X)(0)+X-X(0))(t)\ ,\]
noting that if $X(0)\ge 0$ then $W(X)(0)=X(0)$ and $L(X)(0)=0$. Thus $L_i(X)$ diverges if and only if $L_i(W(X)(0)+X-X(0))$ diverges and, by Theorem~\ref{Theorem:2.1}, if and only if $L_i(X-X(0))$ diverges.  Thus the equivalent conditions in Theorem~\ref{Theorem:2.1} hold for $X$ if and only if they hold for $X-X(0)$. Also, this means that Theorem~\ref{Theorem:2.1} implies that under either of its equivalent conditions, for every $c,b\ge 0$ and even for every $c,b\in\R^n$ we have that
\[
W(c+X)(t)-W(b+X)(t)\r0\ .
\]
\end{Remark}

\begin{Remark}
\rm In the same way that (ii)$\Leftrightarrow$(iii) in Theorem 2.1,
it is tempting to add a fourth equivalence according to which (i)
holds for {\em some} (rather than {\em all}) $a$. This is, however
false. For example, in one dimension with $X(t)=-\min(t,1)$ we have
that $W(a+X)(t)=\max(a-\min(t,1),0)$. Thus, for $a\le 1$ and $t\ge
a$ we have that $W(a+X)(t)=W(X)(t)=0$, but if $a>1$ then
$W(a+X)(t)=a-1>0$ for all $t\ge 1$ while $W(X)(t)=0$, so that (i)
holds for $0\le a\le 1$ but does not hold for $a>1$.
\end{Remark}

It is well known that if $Y$ is nonnegative, nondecreasing right
continuous and $X+RY\ge 0$ then $L(X)\le Y$. That is $L(X)$ is the
minimal nonnegative, nondecreasing right continuous process for
which $X+RL(X)\ge 0$. Denote
\begin{eqnarray}
M_i(X)(t)&=&-\inf_{0\le s\le t}X_i(s)^-\\
N_i(X)(t)&=&-\inf_{0\le s\le t}(R^{-1}X(s))_i^-
\end{eqnarray}
Theorem~4 of \cite{KW1996} is equivalent to $N(X)\le L(X)\le
R^{-1}M(X)$. In addition, since
\begin{equation}
0\le W(X)=X+RL(X)=X+L(X)-P^tL(X)\le X+L(X)
 \end{equation}
 and since $M$ is the minimal nonnegative nondecreasing process for which $X+M(X)\ge 0$, it follows that $M(X)\le L(X)$. Thus we have the following inequality for every $i$ and every $t$:
\begin{equation}
\max(M_i(X),N_i(X))\le L_i(X)\le (R^{-1}M(X))_i
\end{equation}

We thus obtain the following sufficient condition in terms of the
given data $X,P$ of the problem.
\begin{Theorem}\label{Th:3}
If either $M_i(X)$ or $N_i(X)$ is unbounded or equivalently if
\begin{equation}\liminf_{t\r\infty}X_i(t)=-\infty\end{equation} or
\begin{equation}\liminf_{t\r\infty}(R^{-1}X(t))_i=-\infty\end{equation}
(respectively), then for every $a\ge 0$
\begin{equation}
\lim_{t\r\infty}(W_i(a+X)(t)-W_i(X)(t))=0
\end{equation}
\end{Theorem}

In a similar vein, we also have the following necessary condition:
\begin{Theorem}\label{Th:4}
If for every $a\ge 0$
\begin{equation}
\lim_{t\r\infty}(W_i(a+X)(t)-W_i(X)(t))=0
\end{equation}
then $(R^{-1}M(X))_i$ is unbounded and, thus, $M_j(X)$ is unbounded
(that is, $\liminf_{t\r\infty}X_j(t)=-\infty$) for at least one $j$.
\end{Theorem}

It may be noted that, when imposing stochastic assumptions, our asymptotic results above do not require the existence
of a stationary distribution. The following example illustrates this point.

\begin{Example}\label{Ex:BM}
\rm Let $X$ be a standard
$n$-dimensional Brownian motion. Take
$R=I$ (equivalently, $P=0$). It is well known that, almost surely (a.s.),
$\liminf_{t\r\iy}X_i(t)=-\iy$ for
each $1\leq i\leq n$. Hence by Theorem~\ref{Th:3},
$\lim_{t\r\iy } (W(a+X)(t)-W(X)(t))=0$ a.s., for each
$a\geq 0$. It is known that the process $W(X)(\cdot)$ has no
stationary probability distribution in this case; in fact, $W(X)$ is
transient for $n\geq 3$ and is null recurrent for $n=1,2$. As
our conclusion is based on the sufficient condition in \Th{3}, it
seems not to have been noticed earlier for reflected standard Brownian motion with
normal reflection at the boundary. See also \Ex{BM1} .
\end{Example}

The next example shows that our approach does not require Markovian
structure.

\begin{Example}\label{Ex:Renewal}
\rm This example concerns $n$ independent renewal risk processes.
For $i=1,2,\ldots,n$ let \bea X_i(t) & = &   c_it -
\sum_{\ell =1}^{N_i(t)}U^{(i)}_{\ell},~~t\geq 0, \eea
where we make the following hypotheses:

\textbf{(H)} (i) $c_i>0$ is a constant for each $i;$ (ii)
$\{N_i(t):t\geq 0\},\{U^{(k)}_{\ell}:\ell \geq 1 \},1\leq i,k \leq
n$ are independent families of random variables; (iii) for $1\leq
i\leq n,$ $N_i$ is a renewal counting process with positive inter-renewal epochs (thus a.s. having
finitely many jumps on finite time intervals);  let
$A^{(i)}_{\ell},\ell \geq 1$ denote the i.i.d. interarrival times of
$N_i$ with finite positive expectation; (iv) for fixed $k,$
$U^{(k)}_{\ell},\ell \geq 1$ are i.i.d. positive random variables
with finite expectation.

Hence, for $a=(a_1,\ldots,a_n)$, $a_1+X_1,\ldots,a_n+X_n$ are $n$ independent
Sparre Andersen (or renewal risk) processes with respective initial
capitals $a_1,\ldots,a_n$. Clearly, $X_i, 1\leq i\leq n$ are not Markovian in general.

Denote $X=(X_1,\ldots,X_n)$ and take $R=I$ ($P=0$). Then $X_i$ are independent, so neither
$X$ nor $W(a+X)$ are necessarily Markovian.

In addition, suppose for $1\leq i\leq n,$ \bea c_i & = &
\frac{1}{E(A^{(i)}_1)} E(U^{(i)}_1). \eea Then using Theorem 6.3.1
of \cite{RS1999}, we have for each $i$ that, a.s.
$$\limsup_{t\r\iy}X_i(t)=+\iy,~~
\liminf_{t\r\iy}X_i(t)=-\iy\ .$$

Hence by Theorem~\ref{Th:3}, $\lim_{t\r\iy } (W(a+X)(t)-W(X)(t))=0$ a.s., for each $a\geq 0.$ See also \Ex{Re1}.

\end{Example}

\begin{Remark}
\rm It is reasonable to conjecture that the sufficient conditions in
Theorem~\ref{Th:3} are also necessary. Unfortunately, this is not the case.
Let us demonstrate this in a two dimensional setting. More
precisely, consider a feedforward structure with
$P=\left(\begin{array}{cc}0&1\\0&0\end{array}\right)$,
$-X_1(t)=X_2(t)=t|\sin(t)|$. Then $L_1(X)=M_1(X)=N_1(X)$ is
unbounded but $M_2(X)\equiv 0$. Also $N_2(X)\equiv 0$ since
$X_1+X_2\equiv 0$. However
\begin{eqnarray}
L_2(X)(2\pi n)\nonumber &=&-\inf_{0\le s\le 2\pi
n}(X_2(s)-L_1(X)(s))\\ &\ge &L_1(X)(2\pi n)-X(2\pi n)\\ &=&
L_1(X)(2\pi n)\r\infty \nonumber
\end{eqnarray}
as $n\r\infty$ and thus $L_2(X)$ is unbounded. Thus (i) of Theorem
2.1 holds for $i=1,2$ but neither of the two sufficient conditions
of Theorem~\ref{Th:3} holds for $i=2$.
\end{Remark}

\begin{Remark}
\rm Despite the previous remark, it is clear that in the one
dimensional case $L(X)=M(X)=N(X)=R^{-1}M(X)$ and thus for this case
the condition that $\liminf_{t\r\infty}X(t)=-\infty$ is necessary
and sufficient for $W(a+X)(t)-W(X)(t)\r 0$ as $t\r\infty$ for all
$a\ge 0$. This is, of course, obvious for other reasons as well
since the first time that $a+X(t)$ becomes negative is also a time
where the processes $W(a+X)(t)$ and $W(X)(t)$ couple (since both are
zero). From this time and on they remain equal.
\end{Remark}

We conclude this section with the following result concerning
processes which are of interest in queueing networks; see
\cite{K1996,K1997,KW1996,KL2004} and references therein. A process $X$ is said to have stationary increments (in the strong sense) if the law of $\{X(s+t)-X(s)|\ t\ge 0\}$ is independent of $s$.
\begin{Theorem}\label{Th:5}
Let $X$ be a c\`adl\`ag $\mathbb{R}^n$ valued stochastic process having stationary increments. If for its two sided extension (see~\cite{KW1996}) the limit $\xi=\lim_{t\r -\infty}X(t)/t$ exists, is constant and is finite (in particular when
$X$ also has ergodic increments, $E|X(1)-X(0)|<\infty$ and then
$\xi=E(X(1)-X(0))$) and if in addition the stability condition $R^{-1}\xi<0$ holds (coordinate-wise), then

\begin{description}

\item{\rm(i)} $W(X)(t)$ converges in distribution and moreover, for
every a.s. finite random vector $\xi$, $W(\xi+X)(t)$
converges in distribution to the same limit. That is, the limiting
distribution is independent of initial conditions.

\item{\rm(ii)} $W(X)$ has a unique stationary version, in the sense that if, for $\xi$ and $\eta$, $W(\xi+X)$ and
$W(\eta+X)$ are stationary, then the processes $W(\xi+X)$ and $W(\eta+X)$ are identically distributed. For each $t$ the distribution of this stationary version coincides with the
limiting distribution from~(i).

\item{\rm(iii)} Denote $X_s(\cdot)=X(s+\cdot)-X(s)$. There is a uniquely distributed pair $(\xi,X_0)$ ($\xi$ is a random vector and $X_0$ is a process) such that
$(W(\xi+X_0)(s),X_s)$ is distributed like $(\xi,X_0)$ for each $s\ge 0$. Moreover, denoting $W_s(Y)(t)=W(Y)(t+s)$ and $L_s(Y)(t)=L(Y)(s+t)-L(Y)(s)$ we have that the triplet of processes
$(W_s(\xi+X_0),L_s(\xi+X_0),X_s)$
is distributed like $(W(\xi+X_0),L(\xi+X_0),X_0)$. In particular  $W(\xi+X_0)$ is stationary and $L(\xi+X_0)$ has stationary increments.

\item{\rm(iv)} If in addition $J$ is such that with $J_s(\cdot)=J(\cdot+s)$ we have that $(X_s,J_s)\sim (X_0,J)$ (in particular $J$ is stationary and $X$ has stationary increments), then
    \begin{description}
    \item{\rm(i-J)} For every a.s. finite $\xi$, $(W(\xi+X),J)$ converges in distribution to a limit which is independent of $\xi$.
    \item{\rm(ii-J)} $(W(X),J)$ has a unique stationary version.
    \item{\rm(iii-J} There is a uniquely distributed triplet $(\xi,X_0,J)$ such that $(W(\xi+X_0)(s),X_s,J_s)$ is distributed like $(\xi,X_0,J)$ for each $s\ge 0$. Moreover $(W_s(\xi+X_0),L_s(\xi+X_0),X_s,J_s)$ is distributed like $(W(\xi+X_0),L(\xi+X_0),X_0,J)$.
    \end{description}

\end{description}
\end{Theorem}

\proof
\begin{description}
\item{\rm(i)} A direct consequence of \cite{KW1996, K1997}
and Theorem~\ref{Th:3} above.
\item{\rm(ii)} Assume that for some $\xi$ and $\eta$ both
$W(\xi+X)$ and $W(\eta+X)$ are stationary processes. Then, for $0\le
t_1<\ldots<t_m$ and $s>0$, \begin{equation}(W(\xi+X)(t_1),\ldots,W(\xi+X)(t_m))\end{equation} is
distributed like \begin{equation}(W(\xi+X)(t_1+s),\ldots,W(\xi+X)(t_m+s))\end{equation} with a
similar statement when we replace $\xi$ by $\eta$. Clearly,
\begin{equation}(W(\xi+X)(t_1+s),\ldots,W(\xi+X)(t_m+s))\end{equation} converges in distribution
as $s\r\infty$ (as its distribution is independent of $s$) and the
same with $\eta$. According to our results, the difference between
\begin{equation}(W(\xi+X)(t_1+s),\ldots,W(\xi+X)(t_m+s))\end{equation} and
\begin{equation}(W(\eta+X)(t_1+s),\ldots,W(\eta+X)(t_m+s))\end{equation} converges to zero
(coordinate wise, hence in any reasonable norm). Thus these limiting
distributions must be the same and thus
\begin{equation}(W(\xi+X)(t_1),\ldots,W(\xi+X)(t_m))\end{equation} and
\begin{equation}(W(\eta+X)(t_1),\ldots,W(\eta+X)(t_m))\end{equation} are identically
distributed.
\item{(iii)} As was shown in \cite{KW1996}, $X_0$ can be extended to a two sided process, that is, with $t\in\mathbb{R}$ rather than $t\ge 0$. So assume that $X_0$ is already such. It was also shown there that $W(X_{-s})(s)$ is distributed like $W(X)(s)$, is nondecreasing in $s$ and, with the current assumptions, the limit of $W(X_{-s})(s)$ is a.s. finite. Denote this limit by $W^*$. Clearly for every $0\le t_1<\ldots<t_m$ we have that, a.s., \begin{equation}(W(X_{-s})(s),X_0(t_1),\ldots,X_0(t_m))\r (W^*,X_0(t_1),\ldots,X_0(t_m))
    \end{equation}
    and thus, upon shifting forward by $s$ (that is replacing $X_0$ by $X_s$ and $X_{-s}$ by $X_0$) we have by the stationary increments property of $X$ that
    \begin{equation}(W(X_0)(s),X_s(t_1),\ldots,X_s(t_m))\stackrel{d}{\longrightarrow} (W^*,X_0(t_1),\ldots,X_0(t_m))\ .
    \end{equation}
     From our main result it follows that for every a.s. finite $\xi$ we have that $(W(\xi+X_0)(s)-W(X_0)(s))\r 0$ a.s. and thus the difference between $(W(\xi+X_0)(s),X_s(t_1),\ldots,X_s(t_m))$ and $(W(X_0)(s),X_s(t_1)\ldots,X_s(t_m))$ a.s. vanishes as well. Therefore, it also follows that
    \begin{equation}(W(\xi+X_0)(s),X_s(t_1),\ldots,X_s(t_m))\stackrel{d}{\longrightarrow} (W^*,X_0(t_1),\ldots,X_0(t_m))\ .
    \end{equation}
    Hence, if $\xi$ is such that $(W(\xi+X_0)(s),X_s)$ is distributed like $(\xi,X_0)$ then $(W(\xi+X_0)(s),X_s(t_1),\ldots,X_s(t_m))$ is distributed like \\ $(\xi,X_0(t_1),\ldots,X_0(t_m))$, so that necessarily
    \begin{equation}
    (\xi,X_0(t_1),\ldots,X_0(t_m))\sim (W^*,X_0(t_1),\ldots,X_0(t_m))\ ,
    \end{equation}
    as required. Finally, since
    \begin{equation}
    W_s(\xi+X_0)(t)=W(\xi+X_0)(s+t)=W(W(\xi+X_0)(s)+X_s)(t)
    \end{equation}
    and similarly
    \begin{equation}
    L_s(\xi+X_0)(t)=L(W(\xi+X_0)(s)+X_s)(t)
    \end{equation}
    we have that
    \begin{equation}
    \begin{array}{l}
    (W_s(\xi+X_0),L_s(\xi+X_0),X_s)\\ \\ \quad =(W(W(\xi+X_0)(s)+X_s),L(W(\xi+X_0)(s)+X_s),X_s)
    \end{array}
    \end{equation}
    and, as $(W(\xi+X_0)(s),X_s)\sim (\xi,X_0)$, the right hand side is distributed like $(W(\xi+X_0),L(\xi+X_0),X_0)$ as required.
    \item{(iv)} A trivial modification of the proofs of $(i)-(iii)$ and is thus omitted.
    \qed
\end{description}

\begin{Corollary}\label{Cor:Levy}
If $X$ is a multidimensional L\'evy process with $X(0)=0$,
$E|X(1)|<\infty$ and $R^{-1}EX(1)<0$ (coordinate-wise), then the
Markov process $W(X)$ has a (unique) stationary distribution and
converges in distribution to this stationary distribution for every
(a.s. finite) initial condition.
\end{Corollary}

We should be careful to note that we did not show positive Harris
recurrence here, but `only' that the nice properties of Harris
processes hold. Even this, until now was an open question with the
exception of some special cases, such as reflected Brownian motion
and for L\'evy processes for the special case where the inputs are
nondecreasing and are independent, or when the $P$ is of feedforward
type. In the L\'evy case, see, e.g., \cite{K1996, K1997}. Thus an
open problem mentioned in the first paragraph of \cite{KL2004} is
resolved; see also \cite{KW1996}.

We conclude this section with a result involving Markov additive processes with irreducible, positive recurrent countable state Markov modulation. These processes can be intuitively thought of as follows. There is an irreducible, positive recurrent countable state space continuous time
Markov chain $J$ with rate transition matrix $Q$ and stationary distribution vector $\pi$.
When in state $i$, the additive part $X$ behaves like a L\'evy process with characteristic triplet $(c_i,\Sigma_i,\nu_i)$, where $c_i\in\mathbb{R}^n$, $\Sigma_i$ is symmetric positive semidefinite and $\nu_i$ is a (L\'evy) measure with $\int_{\mathbb{R}^n}\min(\|x\|^2,1)\nu_i(dx)<\infty$.
In addition, when $J$ changes states from $i$ to $j$, there are independent jumps (vector valued, possibly a.s. zero) with distribution $G_{ij}$, where $G_{ii}$ is the distribution of the constant vector zero. For a more precise description for the one dimensional case with finite state space modulation, see for example \cite{AK2000}. The multidimensional analogue as well as the countable state space case is defined in an identical manner.

 It is a basic fact that $J$ may be coupled in a.s. finite time $\tau$ with its stationary version. Denoting $X_s(\cdot)=X(s+\cdot)-X(s)$ and $J_s(\cdot)=J(s+\cdot)$, we have that $(X(\tau+t),J(\tau+t))=(X(\tau)+X_\tau(t),J_\tau(t))$ and by the strong Markov property also that
$(X_\tau,J_\tau)$ is a Markov additive process with (jointly) stationary $J_\tau$ and $X_\tau$ having stationary increments.
Note that this implies that if $(W(a+X_\tau)(t)-W(X_\tau)(t))\r 0$ a.s. for every $a$, then also
$(W(X(\tau)+X_\tau)(t),J_\tau(t))-(W(\xi+X_\tau(t)),J_\tau(t))\r 0$ where by (iv) of Theorem~\ref{Th:5} there is a $\xi$ for which $(W(\xi+X_\tau)(t),J_\tau(t))$ is stationary. Thus the following is immediate.

\begin{Corollary}\label{Cor:MAP}
Let $(X,J)$ be a multi-dimensional Markov additive process (see above) where $J$ is an irreducible, positive recurrent countable (or finite) state space continuous time Markov chain with rate transition matrix $Q$ and stationary vector $\pi$. Assume that $\int_{\|x\|>1}\|x\|\nu_i(dx)<\infty$ and $\int_{\mathbb{R}}\|x\|G_{ij}(dx)<\infty$ for all $i$ (finite moment conditions). Denote
$\rho_i=c_i+\int_{\|x\|>1}x\nu_i(dx)\in\mathbb{R}^n$ and $\mu_{ij}=\int_{\mathbb{R}}xG_{ij}(dx)\in \mathbb{R}^n$.
Finally, denote
\begin{equation}
\rho=\sum_i\pi_i\rho_i+\sum_{ij}\pi_iq_{ij}\mu_{ij}
\end{equation}
(so that $X(t)/t\r \rho$ a.s.) and assume that $R^{-1}\rho<0$. Then the Markov process $(W(X),J)$ has a unique stationary distribution and converges to this stationary distribution for every initial condition.
\end{Corollary}

We mention that $\rho_i$ is the rate of the L\'evy process when $J$ is in state $i$. A special case of this is, of course, finite state space Markov modulated (multi-dimensional) reflected Brownian motion with drifts and covariance matrices $(\mu_i,\Sigma_i)$. In this case $\rho=\sum_i\pi_i\mu_i$.

\section{Extensions}
\setcounter{equation}{0}

In this section we consider extensions of Theorem~\ref{Th:3} to some
cases when the drift $b$ and the reflection matrix $R$ need not be
constants. As mentioned in Section 1, the Skorohod problem with
nonconstant drift and reflection has been studied by many authors
in diverse contexts.

We first describe the set-up. Let $(t,\ell,w)\mapsto b(t,\ell,w)$ be
an $\R^n$-valued function and $(t,\ell,w)\mapsto P(t,\ell,w)$ be an
$n\times n$ matrix valued function, where $t$, $\ell\in\R^n$ and
$w\in\R^n$ are all nonnegative. We denote by $b_i$ the $i$th coordinate
of $b$ and by $p_{ij}$ the $i,j$th coordinate of $P$. Finally, let
$R=I-P^t$. In this section, we assume the following conditions.

\begin{description}
\item{(A1)} $b,P$ are bounded, continuous, Lipschitz continuous in
$(\ell,w)$ coordinate-wise uniformly in $t$, $p_{ij}\ge 0$ and
$p_{ii}\equiv 0$.

\item{(A2)} Denoting $\pi_{ij}=\sup_{t,\ell,w}p_{ij}(t,\ell,w)$
(note that $\pi_{ii}=0$) and setting $\Pi=(\pi_{ij})$,
we assume that the spectral radius of $\Pi$ is strictly less than~1.
We term this a {\em uniform spectral radius condition}.

\item{(A3)} $b,R$ are coordinate-wise, nonincreasing in $\ell$ and nondecreasing in
$w$.
\end{description}

Under (A1) and (A2) one can prove (see, e.g., \cite{R2000} and
references therein) that for each c\`adl\`ag $X$, there exists a
unique pair $(W(X),L(X))$ satisfying (\ref{eq:skorohod}) and
conditions S1-S3 (see Section~1).

In this section, the notation $W(a+X),L(a+X)$ becomes a bit
cumbersome and thus, from here on we will abbreviate it using the
notation $W^a=W(a+X)$ and $L^a=L(a+X)$. In particular, $W^0=W(X)$
and $L^0=L(X)$.

If (A1)-(A3) hold, then
$W^a\geq W^0,$ $L^a\leq
L^0,$ and $L^a(t_2)-L^a(t_1)\leq
L^0(t_2)-L^0(t_1)$ for $t_1\leq t_2.$ See Theorems 3.7 and 4.1 of
\cite{R2000}.

\begin{Lemma}
(i) Assume (A1)-(A3). Let $1\leq i\leq n.$ Put \bea \beta_i(s)&=&
\sup\{b_i(s,0,w):w\ge 0\},~~s\geq 0. \eea Suppose \bea
\liminf_{t\rightarrow\infty}[X_i(t)+\int_0^t\beta_i(s)ds] &=&
-\infty. \eea Then for any $a\geq 0$ \bea \lim_{t\rightarrow\infty}
L^a_i(t) &=& +\infty. \eea

(ii) Let (A1)-(A3) hold In addition, let $R$ be a constant
matrix. Let $1\leq i\leq n.$ Put \bea \hat{\beta}_i(s)&=&
\sup\{(R^{-1}b)_i(s,0,w):w\ge 0\},~~s\geq 0. \eea Suppose (3.2) or
\bea \liminf_{t\rightarrow\infty}[(R^{-1}X)_i(t)+\int_0^t
\hat{\beta}_i(s)ds] &=& -\infty, \eea hold. Then (3.3)
holds for any $a\geq 0.$
\end{Lemma}

{\bf Proof:}(i) By (A3) note that $b_i(s,L^a(s),W^a(s))\leq
b_i(s,0,W^a(s))\leq \beta_i(s).$ Therefore, as $p_{ij}\geq
0,i\neq j,$ using (A3), (\ref{eq:skorohod}) it is easy to get \bea
L^a_i(t) &\geq & -a_i-[X_i(t)+\int_0^t\beta_i(s)ds],~~t\geq 0. \eea
As $L^a_i$ is nondecreasing, l.h.s. of (3.3) makes sense; now (3.2),(3.6) imply (3.3).

(ii) Because of part (i) we need to consider only the case when (3.5) holds.
As $R$ is a constant, by (\ref{eq:skorohod}) \Bea L^a(t) &=&
R^{-1}W^a(t) -R^{-1}a -R^{-1}X(t) -\int_0^t
R^{-1}b(s,L^a(s),W^a(s))ds. \Eea Since $R^{-1}$ is
nonnegative, as in part (i) we have \bea L^a_i(t) &\geq &
-(R^{-1}a)_i -[(R^{-1}X)_i(t)+\int_0^t \hat{\beta}_i(s)ds],~~t\geq
0.  \eea Use (3.5),(3.7) now to get (3.3). \qed

We take a closer look at two situations. The first one concerns the
case when the coefficients do not depend on the space variables.

\begin{Proposition}
Let $ b,R$ be functions only of the time variable, satisfying
(A1)-(A3). Then for any $a\geq 0,$ $\lim_{t\rightarrow \infty}
[W^a(t)-W^0(t)]$ exists and belongs to $[0,\iy)^n;$ that is, (2.5)
holds.
\end{Proposition}

{\bf Proof:} Let $\Pi$ be the $n\times n$ constant matrix as in (A2)
with nonnegative entries. Note that $R(s)=I-P^t(s),~s\geq 0;$ also
$\Pi-P^t(s)$ has nonnegative entries. Consequently \bea (I-\Pi)^{-1}R(s)
&=& (I+\Pi+\Pi^2+\cdots)(I-P^t(s)) \nonumber \\ &=&
I+\sum_{k=0}^{\infty}\Pi^k(\Pi-P^t(s)) \eea is a matrix with nonnegative
entries. As the coefficients depend only on the time variable, for
any $a\geq 0,$  \bea && (I-\Pi)^{-1}[W^a(t)-W^0(t)] \nonumber
\\ &=& (I-\Pi)^{-1}a  +
\int_{(0,t]}(I-\Pi)^{-1}R(s)d(L^a-L^0)(s). \eea Because of
(A3), we know that $d(L^a-L^0)\leq 0.$ Hence by
(3.8),(3.9) it now follows that
$(I-\Pi)^{-1}[W^a(t)-W^0(t)]$ is nonincreasing in $t;$ cf.
(iv) of Theorem~\ref{Theorem:1.1} above. As
$W^a-W^0$ is nonnegative, the required
conclusion now follows. \qed

In view of Lemma 3.1, Proposition 3.2 and the proof of
(iii)$\Rightarrow$(i) in Theorem 2.1, the following result is now
immediate.

\begin{Theorem}
Let (A1)-(A3) hold; in addition let $b,R$ be functions
only of the time variable. Suppose \bea \liminf_{t\rightarrow\infty}
[X_i(t)+\int_0^t b_i(s)ds] &=& -\infty,~~1\leq i\leq n. \eea Then $$
\lim_{t\rightarrow\infty} [W^a(t)-W^0(t)] = 0 $$ for any
$a\geq 0.$
\end{Theorem}

{\bf Note:} In addition to the hypotheses of the above theorem,
suppose $R$ is a constant. Then, because of part (ii) of Lemma 3.1,
for the theorem to hold, we need only to assume that for each $i,$
either (3.10) holds or
$$\liminf_{t\rightarrow\infty} [(R^{-1}X)_i(t) +\int_0^t
(R^{-1}b)_i(s)ds] = -\infty. $$

The other situation we consider concerns the special case of
\emph{feedforward} structure. We need some notations. For fixed
$k=1,2,\cdots,n,$ vector $x=(x_1,\cdots,x_n)^t,$ matrix
$A=((A_{ij}))_{1\leq i,j\leq n}$ denote
$x_{[k]}=(x_1,x_2,\cdots,x_k)^t,A_{[k]}=((A_{ij}))_{1\leq i,j\leq
k};$ if $X$ is an $\R^n$-valued function, then
$X_{[k]}=(X_1,\cdots,X_k)^t;$. The following
assumption is characteristic of the feedforward structure in our
context.

{\bf (B1)} For each $1\leq i\leq n,$ $R_{ij}\equiv 0$ if $j > i,$
$R_{ii}\equiv 1;$ $ R_{ij}(s,\ell,w)= R_{ij}(s,\ell_{[i]},w_{[i]})$ if $j
< i,$ $b_i(s,\ell,w) = b_i(s,\ell_{[i]},w_{[i]})$ for $s\geq 0,\ell,w\ge 0.$

So $b_i,R_{i\cdot}$ do not depend on
$\ell_{\ell},w_{\ell},~\ell >i.$ Recalling that $R=I-P^t$, note that
$P_{ij}\equiv 0,j\leq i.$

In this context, note that the so called Skorohod equation (\ref{eq:skorohod}) becomes for
$1\leq i\leq n$ \bea W^a_i(t) &=& a_i +X_i(t) +\int_0^t
b_i(s,L^a_{[i]}(s-),W^a_{[i]}(s-))ds + L^a_i(t)
\nonumber
\\ && +
\sum_{j=1}^{i-1}\int_{(0,t]}R_{ij}(s,L^a_{[i]}(s-),W^a_{[i]}(s-))dL^a_j(s).\eea
Fix $k\in \{1,2,\cdots\}.$ Note that
$L^a_{[k]},W^a_{[k]}$ are functions taking
value in $\R^k$ with $L^a_{[k]}(0)=0,W^a_{[k]}(0)=a_{[k]}.$
By (3.11) it follows that \bea W^a_{[k]}(t) &=& a_{[k]} +
X_{[k]}(t) + \int_0^t
b_{[k]}(s,L^a_{[k]}(s-),W^a_{[k]}(s-))ds \nonumber \\ && +
\int_{(0,t]}R_{[k]}(s,L^a_{[k]}(s-),W^a_{[k]}(s-))dL^a_{[k]}(s).
\eea Clearly $(L^a_{[k]})_i=L^a_i$ is
nondecreasing and can increase only when
$(W^a_{[k]})_i=W^a_i=0,~1\leq i\leq k.$
Summarizing the above and also using uniqueness, we get the
following.

\begin{Proposition} Let (A1)-(A3),(B1) hold. Let $a\ge 0.$ Let
$L^{a },W^a$ be the solution pair for the
Skorohod problem in $[0,\infty)^n$, corresponding to
$(a+X,b,R).$ Fix $1\leq k\leq n.$ Then
$L^a_{[k]},W^a_{[k]}$ is the unique solution
pair for the Skorohod problem in $[0,\infty)^k$ corresponding to
$(a_{[k]}+X_{[k]},b_{[k]},R_{[k]}).$
\end{Proposition}

For our next result we need one more hypothesis.

{\bf (B2)} Coefficients $b,R$ are independent of the
$y$-variables.

\begin{Proposition} Let (A1)-(A3),(B1)-(B2) hold. Let $a\geq 0.$
Assume that (3.3) holds for all $1\leq i\leq n.$ Then there is
$T\geq 0$ such
that \bea L^a(t)-L^a(T) &=& L^0(t)-L^0(T),~~t\geq T,\\
W^a(t) &=& W^0(t), ~~ t\geq T. \eea
\end{Proposition}

{\bf Proof:} Note that it is enough to prove the following: for
fixed $i\in \{1,2,\cdots,n\}$ there exists $t_i\geq 0$ such that
\bea L^a_{[i]}(t)-L^a_{[i]}(t_i) &=&
L^0_{[i]}(t)-L^0_{[i]}(t_i),~~t\geq t_i, \\ W^a_{[i]}(t)
&=& W^0_{[i]}(t),~~t\geq t_i. \eea Then one can take $i=n,T=t_n$
in (3.15),(3.16) to get (3.13),(3.14).

As (3.3) holds for $i=1,$ there is $t\geq 0$ such that $L^a_1(t)
> L^a_1(0)=0.$ Since $L^a_1$ can increase only when
$W^a_1=0,$ there is $t_1\geq 0$ such that
$W^a_1(t_1)=0,$ and hence $W^0_1(t_1)=0$ because of (A3). By
(B1),(B2) clearly $R_{1j}=0,j\geq 2,$
$b_1(s,\ell,w)=\tilde{b}_1(s,w_{[1]})$ where $\tilde{b}_1$ is
an appropriate function. Put
$\tilde{X}_{[1]}(t)=X_{[1]}(t)-X_{[1]}(t_1),$
$\tilde{L}^a_{[1]}(t)=L^a_{[1]}(t)-L^a_{[1]}(t_1),$
$\tilde{L}^{(0)}_{[1]}(t)=L^0_{[1]}(t)-L^0_{[1]}(t_1),t\geq
t_1.$ Hence (\ref{eq:skorohod}) for $i=1$ now implies
$$W^a_{[1]}(t)=0+\tilde{X}_{[1]}(t)+\int_{t_1}^t\tilde{b}_1(s,W^a_{[1]}(s-))ds
+\tilde{L}^a_{[1]}(t),~t\geq t_1,$$
$$W^0_{[1]}(t)=0+\tilde{X}_{[1]}(t)+\int_{t_1}^t\tilde{b}_1(s,W^0_{[1]}(s-))ds
+\tilde{L}^{(0)}_{[1]}(t),~t\geq t_1.$$ Using the above it is not
difficult to see that
$\{(\tilde{L}^a_{[1]}(t),W^a_{[1]}(t)):t\geq t_1\}$ as well
as $\{(\tilde{L}^{(0)}_{[1]}(t),W^0_{[1]}(t)):t\geq t_1\}$ are
both solution pairs to Skorohod problem in $[0,\infty)$ corresponding
to $(0+\tilde{X}_{[1]},b_{[1]},R_{[1]})$ for
$t\geq t_1.$ By uniqueness, we see that (3.15),(3.16) hold for
$i=1.$

We now assume that (3.15),(3.16) hold for $i=1,2,\cdots,k-1$ where
$k\leq n.$ We will now show that they hold for $i=k$ as well.
Without loss of generality we may take $0\leq t_1\leq t_2\leq
\cdots\leq t_{k-1}<\infty.$ Since (3.3) holds for $i=k$ note that
there is $t>t_{k-1}$ such that $L^a_k(t)>L^a_k(t_{k-1})\geq
0.$ So there exists $t_k\geq t_{k-1}$ such that $W^a_k(t_k)=0,$
and hence $W^0_k(t_k)=0.$ As (3.16) holds for $i=k-1,$ we now
have $W^a_{[k]}(t_k)=W^0_{[k]}(t_k)=\zeta_{[k]},$ say. Put
$\tilde{X}_{[k]}(t)=X_{[k]}(t)-X_{[k]}(t_k),$
$\tilde{L}^a_{[k]}(t)=L^a_{[k]}(t)-L^a_{[k]}(t_k),$
$\tilde{L}^{(0)}_{[k]}(t)=L^0_{[k]}(t)-L^0_{[k]}(t_k),t\geq
t_k.$ By Proposition 3.4,
$\{(\tilde{L}^a_{[k]}(t),W^a_{[k]}(t)):t\geq t_k\}$ is the
solution pair for the Skorohod problem in $[0,\infty)^k$ corresponding
to
$(\zeta_{[k]}+\tilde{X}_{[k]},b_{[k]},R_{[k]}).$
By (B2) note that
$b_{[k]}(s,\ell_{[k]},w_{[k]})=\tilde{b}_{[k]}(s,w_{[k]}),$
$R_{[k]}(s,\ell_{[k]},w_{[k]})=\tilde{R}_{[k]}(s,w_{[k]}),$  where
$\tilde{b}_{[k]},\tilde{R}_{[k]}$ are appropriate
functions on $[0,\infty)^{k+1}$. So using (3.12) we get for
$t\geq t_k$  \bea W^a_{[k]}(t) &=& \zeta_{[k]} +
\tilde{X}_{[k]}(t) + \int_{t_k}^t
b_{[k]}(s,\tilde{L}^a_{[k]}(s-)+L^a_{[k]}(t_k),W^a_{[k]}(s-))ds \nonumber \\
&& +
\int_{(t_k,t]}R_{[k]}(s,\tilde{L}^a_{[k]}(s-)+L^a_{[k]}(t_k),W^a_{[k]}(s-))
d\tilde{L}^a_{[k]}(s) \nonumber \\
&=& \zeta_{[k]} + \tilde{X}_{[k]}(t) + \int_{t_k}^t
\tilde{b}_{[k]}(s,W^a_{[k]}(s-))ds \nonumber \\
&&
+\int_{(t_k,t]}\tilde{R}_{[k]}(s,W^a_{[k]}(s-))d\tilde{L}^a_{[k]}(s).
\eea
  In a similar fashion we get for $t\geq t_k$ \bea W^0_{[k]}(t)
  &=& \zeta_{[k]} + \tilde{X}_{[k]}(t) + \int_{t_k}^t
\tilde{b}_{[k]}(s,W^0_{[k]}(s-))ds \nonumber \\
&&
+\int_{(t_k,t]}\tilde{R}_{[k]}(s,W^0_{[k]}(s-))d\tilde{L}^{(0)}_{[k]}(s).
\eea Using (3.17),(3.18) it can be seen that
$\{(\tilde{L}^a_{[k]}(t),W^a_{[k]}(t)):t\geq t_k\}$ as well
as $\{(\tilde{L}^{(0)}_{[k]}(t),W^0_{[k]}(t)):t\geq t_k\}$ are
both solution pairs for Skorohod problem in $[0,\infty)^k$ corresponding
to
$(\zeta_{[k]}+\tilde{X}_{[k]},\tilde{b}_{[k]},\tilde{R}_{[k]})$
for $t\geq t_k;$ the other requirements are very easy to check. By
uniqueness it now follows that (3.15),(3.16) hold for $i=k,$
completing the proof. \qed

Putting together Lemma 3.1 and Proposition 3.5 we get the following
generalization of Theorem 4.1 of \cite{K1997}.

\begin{Theorem} (i) Assume (A1)-(A3),(B1),(B2). Suppose (3.2) holds,
where $\beta_i$ is given by (3.1), for each $1\leq i\leq n.$ Then
for any $a\geq 0$ there is $T\geq 0$ such that
$W^a(t)=W^0(t),~t\geq T.$ In particular
$\lim_{t\rightarrow\infty} [W^a(t)-W^0(t)]=0$ for any $a\geq
0.$

(ii) In addition to the hypotheses in (i), let $R$ be a constant
matrix. Suppose for each $1\leq i\leq n,$ at least one among
(3.2),(3.5) hold. Then the conclusions of part (i) remain valid.
\end{Theorem}

\begin{Example}\label{Ex:BM1}
\rm As in \Ex{BM} let $X$ be a standard
$n$-dimensional Brownian motion and take $R=I.$ Then
by Theorem~3.5, for any $a\geq 0$ there is $T>0$ (depending on $a$
and the sample path) such that $W^a(t)=W^0(t),t\geq T.$ So $W^a$ and
$W^0$ couple in finite time. It may be noted that earlier results
on coupling of reflected Brownian motion in the orthant did not
cover the case of reflected Brownian motion with zero mean and with
normal reflection at the boundary, as the conditions required were
stronger. (This example may be contrasted with the result in
\cite{CL1990}, which says that the reflected Brownian motions in a
planar bounded convex domain (whose curvature is bounded away 0)
with normal reflection, starting from two different points, never
couple in finite time, though the distance between them converges
to zero).
\end{Example}

The next example includes also the case of mean zero and normal
reflection as in \Ex{Renewal}.

\begin{Example}\label{Ex:Re1}
\rm Let $X=(X_1,\ldots,X_n)$ where $X_i,1\leq i\leq n$ are $n$
independent renewal risk processes starting at 0 as in \Ex{Renewal}.
 Let $b,R$
satisfy the hypotheses of Theorem~3.6(i). For $a\geq 0$ let
$W^a(\cdot,\omega)$ denote the regulated/ reflected process in the
orthant corresponding to $a+X$ with coefficients
$b,R$. As $b,R$ do not depend on the
$y$-variables, note that (see Section 6 of \cite{R2000}) $W^a$ is a
strong Markov process. So
by Theorem 3.6(i), for any $a\geq 0,$ there is an almost surely
finite stopping time $T$ such that $W^a(t)=W^0(t)$ for all $t\geq
T;$ in other words $W^a$ and $W^0$ couple in finite time with
probability one. This has the following implication of interest in
actuarial risk theory. Let $\tau^a=\inf\{t>0:W^a(t)=0\}$ denote the
first hitting time of state 0; it is the ruin time (see
\cite{R2011b}) for the multidimensional process $W^a.$ Suppose
$\tau^0 < \iy$ with probability one. Then by the strong Markov property
it follows that $W^0$ visits the state 0 infinitely often with
probability one. Consequently, for any $a\geq 0,$ as $W^a$ and $W^0$
couple in finite time with probability one we have $\tau^a < \iy$
with probability one.
\end{Example}

A crucial ingredient in the proofs of Theorems 2.4, 3.3 is that
\[\lim_{t\rightarrow\infty} [W^a(t)-W^0(t)]\] exists and
belongs to the nonnegative orthant. In Theorem 3.6 this conclusion is a
byproduct of the coupling. So we end with the question: Is there an
analogue of Proposition 3.2 for more general coefficients
$b,R?$

{\bf Acknowledgement} The authors thank two referees and an Associate Editor
for critical comments, helpful suggestions on an earlier draft and for bringing
some relevant references to our notice; these considerably improved the presentation.

\end{document}